\documentclass{amsart}
\usepackage{color}
\usepackage{graphicx,multicol,footmisc}

\newcommand{\bib}[5]{{\bibitem{#1}#2, {\emph{#3},} #4#5.}}

\newcommand{\A}{\mathbb{A}}

\newcommand{\Z}{\mathbb{Z}}
\newcommand{\Q}{\mathbb{Q}}
\newcommand{\R}{\mathbb{R}}
\newcommand{\C}{\mathbb{C}}
\renewcommand{\P}{\mathbb{P}}

\newcommand{\cP}{\mathcal{P}}
\newcommand{\cY}{\mathcal{Y}}
\newcommand{\cM}{\mathcal{M}}
\newcommand{\cH}{\mathcal{H}}

\newcommand{\lra}{\longrightarrow}
\newcommand{\isom}{\approx}
\theoremstyle{plain}
\newtheorem*{probl}{Problem}
\newtheorem{proposition}{Proposition}
\newtheorem{lemma}{Lemma}
\theoremstyle{example}
\newtheorem{example}{Example}

\begin{document}
\parindent0cm

\title{Calabi-Yau conifold expansions}
\author{Slawomir Cynk\and Duco van Straten}
\address{D. van Straten: Institut f\"ur Mathematik, Johannes Gutenberg University,
  55099 Mainz, Germany 
e-mail:\texttt{straten@mathematik.uni-mainz.de}
}
\address{S. Cynk: Institut of Mathematics, Jagiellonian University, 
\L ojasiewicza 6, 30-348 Krak\'ow, Poland,
e-mail:\texttt{slawomir.cynk@uj.edu.pl}
}
\maketitle

\begin{abstract}
We describe examples of computations
of Picard-Fuchs operators for families of Calabi-Yau
manifolds based on the expansion of a period
near a conifold point. We find examples of operators without
a point of maximal unipotent monodromy, thus answering a question posed by J. Rohde.
\end{abstract}

\section{Introduction}
The computation of the instanton numbers $n_d$ for the quintic $X \subset \P^4$
using the period of the quintic mirror $Y$ by P. Candelas, X. de la Ossa and coworkers
\cite{COGP} marked the beginning of intense mathematical interest in the mechanism of mirror symmetry that continues to the present day. 
On a superficial and purely computational level the calculation runs as
follows: one considers the hypergeometric differential operator
\[ \cP=\theta^4-5^5 t (\theta+\tfrac{1}{5})(\theta+\tfrac{2}{5})(\theta+\tfrac{3}{5})(\theta+\tfrac{4}{5})\]
where $\theta=t\frac{d}{dt}$ denotes the logarithmic derivation.
The power series 
\[ \phi(t)=\sum \frac{(5n)!}{(n!)^5}\; t^n\]
is the unique holomorphic solution $\phi(t)=1+\ldots$ to the differential equation
\[ \cP \phi =0\]
There is a unique second solution $\psi$ that contains a log:
\[ \psi(t) =\log(t) \phi(t)+\rho(t)\]
where $\rho \in t \Q[[t]]$. We now define 
\[q:=e^{\psi/\phi}=te^{\rho/\phi}=t+770t^2+\ldots\]
We can use $q$ as a new coordinate and as such it can be used to bring the operator 
$\cP$ into the local normal form
\[\cP=D^2\frac{5}{K(q)}D^2\]
where $D=q\frac{d}{dq}$, and $K(q)$ is a power series.
When we write this series $K(q)$ in the form of a Lambert-series
\[K(q)=5+\sum_{d=1}^{\infty} n_d d^3 \frac{q^d}{1-q^d}\]
one can read off the numbers
\[n_1=2875,\;\; n_2=609250,\;\; n_3=317206375,\ldots\]

The data in the calculation are tied to two Calabi-Yau threefolds:\\

A. The quintic threefold $X \subset \P^4$ ($h^{11}=1, h^{12}=101$). 
The $n_d$ have the interpretation of number of rational degree $d$ curves on
$X$, counted in the Gromov-Witten sense (see \cite{Gi}, \cite{CK}).\\

B. The quintic mirror $Y$ ($h^{11}=101, h^{12}=1$). $Y$ is member
of a pencil $\cY \lra \P^1$, and $\cP$ is Picard-Fuchs operator of
this family. The series $\phi$ is the power-series expansion 
of a special period near the point $0$, which is a point of maximal
unipotent monodromy, a so-called MUM-point.\\

As one can see, the whole calculation depends only on the differential
operator $\cP$ or its holomorphic solution $\phi$, and never uses
any further geometrical properties of $X$ or $Y$, except maybe for
choice of $5$, which is the degree of $X$.\\

In \cite{AZ} this computation was taken as the starting point to 
investigate so-called {\em CY3-operators}, which are fuchsian 
differential operators
$\cP \in \Q(t,\theta)$ of order four with the following properties:

\begin{enumerate}
\item  The operator has the form
\[P=\theta^4+tP_1(\theta)+\ldots+t^rP_r(\theta)\]
 where the $P_i$ are polynomials of degree at most four. 
This implies in particular that $0$ is a MUM-point.\\
\item The operator $\cP$ is {\em symplectic}. This means the $\cP$
leaves invariant a symplectic form in the solution space. The operator
 than is formally self-adjoint, which can be expressed by a simple
condition on the coefficients, \cite{Bo}, \cite{AZ}.\\
\item The holomorphic solution $\phi(t)$ is in $\Z[[t]]$.\\ 
\item Further integrality properties: the expansion of the
$q$-coordinate has integral coefficients and the instanton numbers
are integral (possibly up to a common denominator).
\end{enumerate}

There is an ever growing list of operators satisfying the first three 
and probably the last conditions \cite{AESZ}. It starts with the
above operator and continues with $13$ further hypergeometric
cases, which are related to Calabi-Yau threefolds that are complete intersections
in weighted projective spaces. Recently, M. Bogner and S. Reiter \cite{Bo}, \cite{BR} have classified and constructed the symplectically rigid Calabi-Yau operators, thus
providing a solid understanding for the beginning of the list.

Another nice example is operator nr. $25$ from the list:
\[\cP=\theta^4-4t(2\theta+1)^2(11\theta^2+11\theta+3)-16t^2(2\theta+1)^2(2\theta+3)^2\]
The holomorphic solution of the operator is  $\phi(t)=\sum A_nt^n$ 
where
\[A_n:={2n \choose n}^2 \sum_{k=1}^n {n \choose k}^2{n+k \choose k}\]
This operator was obtained in \cite{BCKS} as follows: one considers
the Grassmanian $Z:=G(2,5)$, a Fano manifold of dimension $6$, with
$Pic(Z) \approx \Z$, with ample generator $h$, the class of a hyperplane section in
the Pl\"ucker embedding. As the canonical class of $Z$ is $-5h$, the
complete intersection $X:=X(1,2,2)$ by hypersurfaces of degree $1,2,2$ is
a Calabi-Yau threefold with $h^{11}=1, h^{12}=61$. The small quantum-cohomology 
of $Z$ is known, so that one can compute its quantum D-module. The quantum Lefschetz
theorem then produces the above operator nr. $25$ which thus provides the
numbers $n_d$ for $X$:
\[ n_1=400,\;\; n_2 = 5540,\;\; n_3 = 164400,\;\ldots\]
 Also, a mirror manifold $Y=Y_t$ was described as
(the resolution of the toric closure of) a hypersurface in the torus $(\C^*)^4$
given by a Laurent-polynomial.\\

The question arises which operators in the list are related in a similar
way to a mirror pair $(X,Y)$ of Calabi-Yau threefolds with $h^{11}(X)=h^{12}(Y)=1$. 
This is certainly not to be expected for all operators, but it suggests the following attractive problem.\\
\begin{probl}
A. Construct examples of Calabi-Yau threefolds $X$ with $h^{11}=1$ and try
to identify the associated quantum differential equation.\\

B. Construct examples of pencils of Calabi-Yau threefolds $\cY \lra \P^1$ 
with $h^{12}(Y_t)=1$ and try to compute the the associated Picard-Fuchs 
equation.\\

It has been shown that in many cases one can predict from the operator $\cP$ 
alone topological invariants of $X$ like $(h^3, c_2(X)h, c_3(X)$), \cite{ES} 
and the zeta-function of $Y_t$ \cite{SS}, \cite{Yu}. 
In either case we see that the operators of the list provide
predictions for the existence of Calabi-Yau threefolds with quite 
precise properties. Recently, A. Kanazawa \cite{Ka} has used weighted
Pfaffians to construct some Calabi-Yau threefolds $X$ whose existence were
predicted in \cite{ES}. In this note we report on work in progress
to compute the Picard-Fuchs equation for a large number of families
of Calabi-Yau  threefolds with $h^{12}=1$.
\end{probl}

\section{How to compute Picard-Fuchs operators}

{\bf The method of Griffiths-Dwork}\medskip

For a smooth hypersurface $Y \subset \P^{n}$ defined by
a polynomial $F \in \C[x_0,\ldots,x_n]$ of degree $d$ one 
has a useful representation of (the primitive part of) the
middle cohomology $H^{n-1}_{prim}(Y)$ using residues of differential
forms on the complement $U:=\P^n\setminus Y$. One can work with the
complex of differentials forms with poles along $Y$ and compute
modulo exact forms. Although this method was used in the $19$th
century by mathematicians like Picard and Poincar\'e, it was first
developed in full generality by P. Griffiths \cite{Gr} and B. Dwork \cite{Dw}
in the sixties of the last century.\\
The {\em Griffths isomorphism} identifies the Hodge space $H^{p,q}_{prim}$ 
with a graded piece of the Jacobian algebra 
\[R:=\C[x_0,\ldots,x_n]/(\partial_0F,\partial_1F,\ldots, \partial_nF)\]
More precisely one has
\[
\begin{array}{rcl}
 R_{d(k+1)-(n+1)} &\stackrel{\isom}{\lra}& H^{n-1-k,k}_{prim}(Y)\\[2mm]
 P &\mapsto & Res(\frac{P\Omega}{F^{k+1}})
\end{array}
\]
where $\Omega:=\iota_E (dx_0\wedge dx_1 \wedge\ldots dx_n)$,
and $E=\sum x_i\partial/\partial x_i$ is the Euler vector field.
This enables us to find an explicit basis.\\

If the polynomial $F$ depends on a paramter $t$, we obtain a
pencil $\cY \lra \P^1$ of hypersurfaces, which can be seen
as a smooth hypersurface $Y_t$ over the function 
field $K:=\C(t)$ and the above method provides a basis 
$\omega_1,\ldots,\omega_r$ of differential
forms over $K$. We now can differentiate the differential forms
 $\omega_i$ with respect to $t$ and express the result in
the  basis. This step involves a Gr\"obner-basis calculation.
As a result we obtain an $r\times r$ matrix $A(t)$ with entries
in $K$ such that
\[\frac{d}{dt} \left( \begin{array}{c} \omega_1\\\omega_2\\\ldots\\\omega_r\end{array}
\right)= A(t) \left( \begin{array}{c} \omega_1\\\omega_2\\\ldots\\\omega_r
\end{array}\right)
\]
The choice of a {\em cyclic vector} for this differential system then provides
a differential operator $\cP \in \C(t,\theta)$ that annihilates all 
period integrals $\int_{\gamma} \omega$. In the situation of Calabi-Yau manifolds
there is always a natural vector obtained from the holomorphic differential.
For details we refer to the literature, for example \cite{CK}.

This methods works very well in simple examples and has been used
by many authors. It can be generalised to the case of (quasi)-smooth 
hypersurfaces in weighted projective spaces and more generally complete 
intersections in toric varieties \cite{BC}. Also, it is possible to handle
families depending on more than one parameter. A closely related method for
tame polynomials in affine space has been implemented by M. Schulze
\cite{Sc1} and H. Movasati \cite{Mo} in {\sc Singular}. The ultimate
generalisation of the method would be an implementation of the
{\em direct image functor in the category of $D$-modules}, which in principle
can be achieved by Gr\"obner-basis calculations in the Weyl-algebra.\\

The {\sc Griffiths-Dwork} method however also has some drawbacks:

\begin{itemize}
\item In many situation the varieties one is interested in 
have singularities. For the simplest types of singularities it
is still possible to adapt the method to take the singularities 
into account, but the procedure becomes increasingly cumbersome
for more complicated singularities.\\
\item In many situations the variety under consideration
is given by some geometrical construction and a description 
with equations seems less appropriate.\\
\end{itemize}
   
In some important situations the following {\em alternative method}
can be used with great succes.\\

{\bf Method of Period Expansion}\medskip

In order to find Picard-Fuchs operator for a family $\cY \lra \P^1$
one does the following:

\begin{itemize}
\item Find the explicit power series expansion of a single period
\[\phi(t)=\int_{\gamma_t}\omega_t =\sum_{n=0}^{\infty} A_n t^n\]

\item Find a differential operator 
\[\cP = P_0(\theta)+tP_1(\theta)+\ldots +t^rP_r(\theta)\]
that annihilates $\phi$ by solving the linear resursion
\[ \sum_{i=0}^r P_i(n)A_{n-i}=0\]
on the coefficients. Here the $P_i$ are polynomials in $\theta$ of 
a certain degree $d$. As $\cP$ contains $(d+1)(r+1)$ coefficients, 
we need the expansion of $\phi$ only up to sufficiently high order 
to find it. 
\end{itemize}

This quick-and-dirty method surely is very old and goed back to the
time of {\sc Euler}. And of course, many important issues arise like: 
{\em To what order do we need to compute our period?} For this one needs 
a priori estimates for $d$ and $r$, which might not be
available. Or: {\em Is the operator $\cP$ really the Picard-Fuchs operator of the family ?} 
We will not discuss these issues here in detail, as they are not so important
in practice: one expands until one finds an operator and if the monodromy
representation is irreducible, the operator obtained is necessarily the 
Picard-Fuchs operator.\\

However, it is obvious that the method stands or falls with our ability to 
find such an explicit period expansion. It appears that the 
{\em critical points} of our family provide the clue.\\

\medskip
\centerline{\bf \em Principle}
\smallskip

{\em If one can identify explicitely a vansihing cycle, then its period
can be computed ``algebraically''.}\\

If our family $\cY \lra \P^1$ is defined over $\Q$, or more generally
over a number field, then it is known that such expansions are $G$-functions
and thus have very strong arithmetical properties, \cite{An}.

Rather then trying to prove here a general statement in this direction,
we will illustrate the principle in two simple examples. The appendix
contains a general statement that covers the case of a variety 
aquiring an ordinary double point.\\

{\bf I.} Let us look at the {\sc Legendre} family of elliptic curves
given by the equation

\[ y^2=x(t-x)(1-x)\]

If the parameter $t$ is a small positive real number, the
real curve contains a cycle $\gamma_t$ that runs from $0$ to $t$ 
and back. If we let $t$ go to zero, this loop shrinks to a point and the
curve aquires an $A_1$ singularity.
The period of the holomorphic differential $\omega=dx/y$ along this
loop is 
\[ \phi(t) =\int_{\gamma_t} \omega=2 F(t)\]
where
\[F(t):=\int_0^t\frac{dx}{\sqrt{(x(t-x)(1-x)}}\]
By the substitution $x \mapsto tx$ we get

\[F(t)=\int_0^1\frac{1}{\sqrt{(1-xt)}}\frac{dx}{\sqrt{x(1-x)}}\]

The first square root expands as 
\[\frac{1}{\sqrt{(1-xt)}}=\sum_{n=0}^{\infty} {2n \choose n}\left(\frac{xt}{4}\right)^n\]
so that
\[F(t)=\sum_{n=0}^{\infty} {2n \choose n} \left(\int_0^{1}
\frac{x^n}{\sqrt{x(1-x)}}dx\right) \,t^{n}\]
The appearing integral is well-known since the work of {\sc Wallis} and is a 
special case of {\sc Euler}s Beta-integral. 
\[ \int_0^1 \frac{x^n}{\sqrt{(x(1-x)}}dx=\pi {2n \choose n} \frac{1}{4^n}.\]
So  the final result is the beautiful series
\begin{eqnarray*}
  F(t)&=&\pi \sum_{n=0}^{\infty} {2n \choose n}^2 {\left(\frac{t}{16}\right)}^n\\[1mm]
&=&\pi{\left(1+{\left(\frac{1}{2}\right)}^2t+{\left(\frac{1\cdot 3}{2\cdot
              4}\right)}^2t^2 + {\left(\frac{1\cdot 3 \cdot 5}{2 \cdot 4 \cdot 6}\right)}^2t^3+\ldots\right)}
\end{eqnarray*}
From this series it is easy to see that the second order operator with 
$F(t)$ as solution is:
\[ 4 \theta^2-t(2\theta+1)^2\]
In fact, the first six coefficients suffice to find the operator.

This should be compared to the {\sc Griffiths-Dwork} method, which would consist
of considering the basis 
\[\omega_1=dx/y,\omega_2=xdx/y\]
of differential forms on $E_t$ and expressing the derivative
\[ \partial_t \omega_1 = - \frac{x(1-x) dx}{(x(t-x)(1-x))^{3/2}}\]
in terms of $\omega_1, \omega_2$ modulo exact forms.\\

{\bf II.} In mirror symmetry one often encounters families of Calabi-Yau manifolds 
that arise from a Laurent polynomial 
\[ f \in \Z[x_1,x_1^{-1},x_2,x_2^{-1},\ldots,x_n,x_n^{-1}]\]
Such a Laurent polynomial $f$ determines a family of hypersurfaces in a torus
given by
\[ V_t:=\{1-tf(x_1,\ldots,x_n)=0\} \subset (\C^*)^{n}\]
In case the Newton-polyhedron $N(f)$ of $f$ is {\em reflexive}, a crepant 
resolution  of the closure of $V_t$ in the toric manifold
determined by $N(f)$ will be a Calabi-Yau manifold $Y_t$. To compute its
Picard-Fuchs operator, the {\sc Griffiths--Dwork} method is usually cumbersome.\\

The holomorphic $n-1$-form on $Y_t$ is given on $V_t$ 
\[ \omega_t:=Res_{V_t}\left(\frac{1}{1-tf}
  \frac{dx_1}{x_1}\frac{dx_2}{x_2} \ldots \frac{dx_n}{x_n}\right)\] 
There is an $n-1$-cycle $\gamma_t$ on $V_t$ whose {\sc Leray}-coboundary is homologous
to $T:=T_{\epsilon}:=\{|x_i|=\epsilon\} \subset (\C^*)^n$. The so-called {\em principal period} is 

\begin{eqnarray*}
\phi(t)&=&\int_{\gamma_t} \omega_t
=\frac{1}{(2\pi i)^n}\int_T \frac{1}{1-tf}
\frac{dx_1}{x_1}\frac{dx_2}{x_2} \ldots \frac{dx_n}{x_n} 
=\sum_{n=0}^{\infty}[f^n]_0 t^n
\end{eqnarray*}

where $[g]_0$ denotes the {\em constant term} of the Laurent series $g$. For this 
reason, the series $\phi(t)$ is sometimes called the  {\em constant term series} 
of the Laurent-polynomial. This method was used in \cite{BS} to determine the 
Picard-Fuchs operator for certain families $Y_t$ and has been popular ever since. 
A fast implementation for the computation of $[g]_0$ was realised by P. Metelitsyn, \cite{Mt}.

\section{Double Octics}
One of the simplest types of Calabi-Yau threefolds is the so-called
{\em double octic}, which is a double cover $Y$ of $\P^3$ ramified over
a surface of degree $8$. It can be given by an equation of the form
\[ u^2=f_8(x,y,z,w)\]
and thus can be seen as a hypersurface in weighted projective space
$\P(1^{4},4)$. For a general choice of $f_8$ the variety $Y$ is 
smooth and has Hodge numbers $h^{11}=1,h^{12}=149$. A nice sub-class of
such double octics consists of those for which  $f_8$ is a product 
of eight planes. In that case $Y$ has singularities at the intersections 
of the planes. In the generic such situation $Y$ is singular along 
$8.7/2=28$ lines, and by blowing up these lines (in any order) we 
obtain a smooth Calabi-Yau manifold $\tilde{Y}$ with $h^{11}=29, h^{12}=9$. 
By taking the eight planes in special positions, the double cover $Y$ aquires other singularities
and a myriad of different Calabi-Yau threefolds with various Hodge numbers
appear as crepant resolutions $\tilde{Y}$. In \cite{Me} $11$ configurations leading to rigid Calabi-Yau
varities were identified. Furthermore, C. Meyer listed $63$ 1-parameter
families which thus give $63$ special 1-parameter families of Calabi-
Yau threefolds $\tilde{Y}_t$, and it is for these that we want to
compute the associated Picard-Fuchs equation. Due to the singularities
of $f_8$, a Griffiths-Dwork approach is cumbersome, if not impossible.
So we resort to the {\em period expansion method}.\\
 
In many of the $63$ cases one can identify a {\em vanishing tetrahedron}: for
a special value of the parameter one of the eight planes passes through
a triple point of intersection, caused by three other planes. In
approprate coordinates we can write our affine equation as
\[ u^2=xyz(t-x-y-z)P_t(x,y,z)\]   
where $P_t$ is the product of the other four planes and we assume
$P_0(0,0,0) \neq 0$.
Analogous to the above calculation with the elliptic curve we now ``see''
a cycle $\gamma_t$, which consists of two copies of the real
tetrahedron $T_t$ bounded by the plane $x=0$, $y=0$, $z=0$, $x+y+z=t$. For
$t=0$ the terahedron shrinks to a point. So we have
\[ \phi(t)=\int_{\gamma_t} \omega=2F(t)\]
where
\[F(t)=\int_{T_t} \frac{dxdydz}{\sqrt{(xyz(t-x-y-z)P_t(x,y,z)}}\]

\begin{proposition}
The period $\phi(t)$ expands
in a series of the form
\[ \phi(t) =\pi^2 t (A_0+A_1t+A_2t^2+\ldots)\]
with $A_i \in \Q$ if $P_t(x,y,z) \in \Q[x,y,z,t]$, $P_{0}(0,0,0)\neq0$.\\
  
\end{proposition}

{\bf proof:} When we replace  $x,y,z$ by $tx,ty,tz$ respectively, we obtain an
integral over the standard tetrahedron $T:=T_1$:
\[F(t)=t\int_T \frac{dxdydz}{\sqrt{xyz(1-x-y-z)}} \frac{1}{\sqrt{P_t(tx,ty,tz)}}\]
 
We can expand the last square root in a power series
\[\frac{1}{\sqrt{P_t(tx,ty,tz)}}=\sum_{iklm}C_{iklm}x^ky^lz^mt^i\]\
and thus find $F(t)$ as a series
\[F(t)=t \sum_{i,k,l,m} \int_T \frac{x^ky^lz^m dx dy dz}{\sqrt{xyz(1-x-y-z)}} C_{iklm} t^i \]

The integrals appearing in this sum can be evaluated easily in terms of the\\
{\bf Generalised Beta-Integral}
\[\int_T x_1^{\alpha_1-1}x_2^{\alpha_2-1}\ldots x_n^{\alpha_n-1}(1-x_1-\ldots-x_n)^{\alpha_{n+1}-1}dx_1dx_2\ldots dx_n =\]
\[ \Gamma(\alpha_1)\Gamma(\alpha_2) \ldots \Gamma(\alpha_{n+1})/\Gamma(\alpha_1+\alpha_2+\ldots+\alpha_{n+1})\]

In particuler we get

\begin{eqnarray*}
\int_T \frac{x^ky^lz^m dx dy dz}{\sqrt{xyz(1-x-y-z)}}&=&\frac{\Gamma(k+1/2)\Gamma(l+1/2) \Gamma(m+1/2)\Gamma(1/2)}{\Gamma(k+l+m+2)}\\
&=&\pi^2 \frac{(2k)!(2l)!(2m)!}{4^{k+l+m} k!l!m!(k+l+m+1)!} \in \pi^2\Q
\end{eqnarray*}
and thus we get an expansion of the form
\[F(t)=\pi^2 t (A_0+A_1t+A_2t^2+A_3t^3+\ldots)\]
where $A_i \in \Q$ when $P_t(x,y,z) \in \Q[x,y,z,t]$ \hfill $\diamond$\\

\begin{example}
Configuration no. $36$ of C. Meyer (\cite{Me}, p.57)
is equivalent to the double octic with equation
\[ u^2=xyz(t-x-y-z)(1-x)(1-z)(1-x-y)(1+(t-2)x-y-z)\]
A smooth model has $h^{11}=49, h^{12}=1$. For $t=0$ the resolution is
a rigid Calabi-Yau with $h^{11}=50, h^{12}=0$, corresponding to arrangement
no. $32$. 
The expansion of the tetrahedral integral around $t=0$ reads:
\[ F(t)=\pi^2 t (1+t+{\frac {43}{48}}\,{t}^{2}+{\frac {19}{24}}\,{t}^{3}+{\frac {10811}
{15360}}\,{t}^{4}+{\frac {9713}{15360}}\,{t}^{5}+\ldots)\]
The operator is determined by the first $34$ terms of the expansion and
reads
\[32\,\theta \left( \theta-2 \right)  \left( \theta-1 \right) ^{2}-16\,t \theta \left( \theta-1 \right)  \left( 9\,{\theta}^{2}-13\,\theta+8 \right)+\]
\[8\,t^2 \theta \left( 33\,{\theta}^{3}-32\,{\theta}^{2}+38\,\theta-10 \right)-t^3(252\,{\theta}^{4}+104\,{\theta}^{3}+304\,{\theta}^{2}+76\,\theta+20)+\]
\[t^4(132\,{\theta}^{4}+224\,{\theta}^{3}+292\,{\theta}^{2}+160\,\theta+38)-t^5(36\,{\theta}^{4}+104\,{\theta}^{3}+140\,{\theta}^{2}+88\,\theta+21)+\]
\[4\, t^6 \left( \theta+1 \right) ^{4}\]

The Riemann symbol of this operator is
\[\left\{ 
\begin{array}{cccc} 
0&1&2& \infty\\
\hline
0&0&0&1\\
1&0&0&1\\
1&0&2&1\\
2&0&2&1\\
\end{array}
\right\}\]
At $0$ we have indeed a 'conifold point' with its characteristic
exponents $0,1,1,2$. At $t=1$ and $t=\infty$ we find MUM-points.
M. Bogner has shown that via a quadratic transformation this operator 
can be transformed to operator number $10^*$ from the AESZ-list, which has
Riemann-symbol
\[\left\{ 
\begin{array}{ccc} 
0&1/256& \infty\\
\hline
0&0&1/2\\
0&0&1\\
0&1&1\\
0&1&3/2\\
\end{array}
\right\}\]
which is symplectically rigid, \cite{BR}. So the family of double 
octics provides a clean $B$-interpretation for this operator.\\
  
\end{example}
\begin{example}
Configuration no. $70$ of Meyer is isomorphic to
\[u^2=xyz(x+y+z-t)(1-x)(1-z)(x+y+z-1)(x/2+y/2+z/2-1)\]
Again, for general $t$ we obtain a Calabi-Yau 3-fold with $h^{11}=49, h^{12}=1$
and for $t=0$ we have $h^{11}=50, h^{12}=0$, corresponding to the rigid Calabi-Yau
of configuration no. $69$ of \cite{Me}. The tetrahedral integral expands as 
\[F(t)=\pi^2 t(1+{\frac {13}{16}}\,t+{\frac {485}{768}}\,{t}^{2}+{\frac {12299}{24576
}}\,{t}^{3}+{\frac {534433}{1310720}}\,{t}^{4}+{\frac {21458473}{
62914560}}\,{t}^{5}+\ldots)\]
and is annihilated by the operator
\[16\,\theta \left( \theta-2 \right)  \left( \theta-1 \right) ^{2}-2\,t \theta \left( \theta-1 \right)  \left( 24\,{\theta}^{2}-24\,\theta+13 \right) +\]
\[t^2 {\theta}^{2} \left( 52\,{\theta}^{2}+25 \right)-2\, t^3 \left( 3\,{\theta}^{2}+3\,\theta+2 \right)  \left( 2\,\theta+1 \right) ^{2} +\] 
\[t^4 \left( 2\,\theta+1 \right)  \left( \theta+1 \right) ^{2} \left( 2\,\theta+3 \right) \]
The Riemann symbol of this operator is:
\[\left\{ 
\begin{array}{cccc} 
0&1&2& \infty\\
\hline
0&0&0&1/2\\
1&0&0&1\\
1&1&1&1\\
2&1&1&3/2\\
\end{array}
\right\}\]
so we see that it has {\em no point of maximal unipotent monodromy}!\\

The first examples of families Calabi-Yau manifolds without MUM-point
were described by J. Rohde \cite{Ro} and studied further by A. Garbagnati
and B. van Geemen \cite{GG}. It should be pointed out that in those cases the 
associated Picard-Fuchs operator was of second order, contrary to the
above fourth order operator. M. Bogner has checked that this operator 
has $Sp_4(\C)$ as differential Galois group. It is probably one of the 
simplest examples of this sort. J. Hofmann has calculated with his
package \cite{Ho} the integral monodromy of the operator. 
In an appropriate basis it reads

\begin{eqnarray*}
T_2 = && \left(\begin{array}{rrrr}\tabcolsep=2cm
    2& \;-7& \phantom{-}\;1& \phantom{-}0\\  \phantom{-}0& 1& 0& 1\\  -1&
    7& 0& 7\\  0& 0& 0& 1\end{array}\right), \ 
\;T_1=  \left(\begin{array}{rrrr}  -1&\,\,\,\, -2& \;\phantom{-}0&
   \; \phantom{-}0\\  2& 3& 0& 0\\  11& 7& 2& 1\\ 
  -3& 1& -1& 0\end{array}\right), \
\\[1mm]T_0=&&  \left(\begin{array}{rrrr}  \phantom{-}1& \phantom{-}0&
  \phantom{-}0& \phantom{-}0\\  0& 1& 0& 0\\  0& 0& 1& 0\\  8& 
  -16& 4& 1\end{array}\right), \  
T_{\infty} = \left(\begin{array}{rrrr}  \phantom{-}0&\;\phantom{-} 23&
    \; -3&\; -2\\  0&
    -15& 2& 1\\  0& -84& 11& 6\\  1& -75& 11& 4 \end{array}\right)   
\end{eqnarray*}

with $T_2 T_1 T_0 T_\infty = id$.\\
\end{example}

As Calabi--Yau operators in the sense of \cite{AZ} need a to have a MUM, 
W. Zudilin has suggested to call an operator without such a point of
maximal unipotent monodromy an {\em orphan}.\\
 
\begin{example}
Configuration no. $254$ of C. Meyer gives a family of 
Calabi-Yau threefolds with $h^{11}=37, h^{12}=1$:
\[u^2=xyz(t-x-y-z)P_t(x,y,z)\]
with
\[
\begin{array}{rcl}
P_t(x,y,z)&=&(1-3z+t-t^2x+tz-tx-2y)(1-z+tx-2x) \cdot\\
        & &\cdot (1-tx+z)(1+t-t^2x+tz-5tx+z-2y-4x)
\end{array}
\]
For $t=0$ we obtain the rigid configuration no. $241$ with $h^{11}=40, h^{12}=0$. The
tetrahedral integral expands as
\[ F(t)=\pi^2 t (
1+\frac{1}{2}\,t+{\frac {37}{24}}\,{t}^{2}+{\frac {41}{16}}\,{t}^{3}+{\frac {
13477}{1920}}\,{t}^{4}+{\frac {14597}{768}}\,{t}^{5}+\ldots)\]

The operator is very complicated and has the following Riemann symbol:
\[\left\{ 
\begin{array}{ccccccccc} 
\alpha_1&\alpha_2&0&\rho_1&\rho_2&\rho_3&-1&1& \infty\\
\hline
0&0&0&0&0&0&0&0&3/2\\
1&1&1&1&1&1&0&0&3/2\\
1&1&1&3&3&3&0&0&3/2\\
2&2&2&4&4&4&0&0&3/2\\
\end{array}
\right\}\]
where at $0$ and $\alpha_{1,2}=-2\pm\sqrt{5}$ we find conifold points,
at the $\rho_{1,2,3}$, roots of the cubic equation $2t^3-t^2-3t+4=0$
we have apparent singularities and at $-1,1$ we find point of
maximal unipotent monodromy, which we also find at $\infty$, after
taking a square root. This operator was not known before.\\
\end{example}

These three examples illustrate the current {\em win-win-win} aspect
of these calculations. It can happen that the operator is known, in which
case we get a nice geometric incarnation of the differential equation.
It can happen that the operator does not have a MUM-point, in which case
we have found a further example of of family of Calabi-Yau threefolds without
a MUM-point. From the point of mirror-symmetry these cases are of special
importance, as the torus for the SYZ-fibration, which in the ordinary cases 
vanishes at the MUM-point, is not in sight. Or it can happen that we find
a {\em new} operator with a MUM-point, thus extending the AESZ-table \cite{AESZ}.  
  
Many more examples have been computed, in particular also for
other types of families, like fibre products of rational elliptic
surfaces of the type considered by C. Schoen, \cite{Sc2}. The first
example of $Sp_4(\C)$-operators without MUM-point were found among these,
\cite{St}. A paper collecting our results on periods of double octics 
and fibre products is in preparation, \cite{CS}.

\section{An algorithm}

Let $\cY$ be a smooth variety of dimension $n$ and $f:\cY \lra \P^1$ a non-constant map to $\P^1$ and let $P \in \cY$ be a critical point. 
In order to analyze the local behaviour of periods of cycles vanishing at $P$, 
we replace $\cY$ by an affine part, on which we have a function $f:\cY \lra \A^1$, with 
$f(P)=0$. An $n$-form \[ \omega \in \Omega^n_{\cY,P} \]
gives rise to a family of differential forms on the fibres of $f$:
\[ \omega_t:=Res_{Y_t}(\frac{\omega}{f-t})\]\\
The period integrals
\[ \int_{\gamma_t} \omega_t\]
over cycles $\gamma_t$, vanishing at $P$ only depend on the class
of $\omega$ in the {\em Brieskorn module} at $P$, which is defined as
\[ \cH_P:=\Omega_{\cY,P}^n/df\wedge d\Omega_{\cY,P}^{n-2}\]
If $P$ is an isolated critical point, it was shown in \cite{Br} that
the completion $\widehat{\cH_P}$ is a (free) $\C[[t]]$-module of rank 
$\mu(f,P)$, the Milnor number of $f$ at $P$. In particular, if $f$ has
an $A_1$-singularity at $P$, we have $\mu (f,P)=1$, and the image of
the class of $\omega$ under the isomorphism
$\widehat{\cH_P} \lra \C[[t]]$ is, up to a factor, just the expansion of
the integral of the vanishing cycle. We will now show how one can calculate
this with a simple algorithm.\\

\begin{proposition}
If $f: \cY \lra \A^1$ and the critical point $P$ of type $A_1$.
If $f:\cY \lra \A^1$, $P$  and $\omega \in \Omega_{\cY,L}$ are defined over $\Q$,
then the period integral over the vanishing cycle $\gamma(t)$
\[\phi(t)=\int_{\gamma(t)}\omega_t\]  
has an expansion of the form
\[\phi(t) = c t^{n/2-1}(1+A_1t+A_2t^2+\ldots)\]
where 
\[ c=d \frac{n}{2} \frac{\Gamma(1/2)^n}{\Gamma(\frac{n}{2}+1)} \]
where $d^2 \in \Q$ and the $A_i\in \Q$ can be computed via a simple algorithm.\\
\end{proposition}

{\bf proof:} As $P$ and $f$ are defined over $\Q$, we may assume that in appropriate
formal coordinates $x_i$ on $\cY$ we have $P=0$, $f(P)=0$ and the map is represented
by a series
\[f=f_2+f_3+f_4+\ldots \]
where $f_2$ is a non-degenerate quadratic form and the $f_d \in \Q[x_1,\ldots,x_n]$ are homogeneous polynomials of degree $d$. After a linear coordinate transformation (which may involve a quadratic field extension) we may and will assume that
\[f_2=x_1^2+x_2^2+\ldots+x_n^2\] 
For $t >0$ small enough, the part of solution set $\{(x_1,x_2,\ldots,x_n) \in \R^n\;|\;f=t\}$  
near $0$ looks like a slightly bumped sphere $\gamma(t)$ and is close to standard sphere
$\{(x_1,x_2,\ldots,x_n) \in \R^n\;|\;f_2=t\}$. This is the vanishing cycle we want to
integrate $\omega_t=Res(\Omega/(f-t))$ over. Note that
\[\int_0^t \int_{\gamma(t)}\omega_t=\int_{\Gamma(t)} \omega\]
where 
\[\Gamma(t)=\cup_{s \in [0,t]} \gamma(s)= \{(x_1,x_2,\ldots,x_n) \in \R^n\;|\;f \le t \}\]
is the {\em Lefschetz thimble}, which is a slightly bumped ball, that is near to the {\em standard ball}
\[B(t):=\{(x_1,x_2,\ldots,x_n) \in \R^n\;|\;f_2 \le t\}\]
The idea is now to change to coordinates that map $f$ into its quadratic part
$f_2$. An automorphism $\phi:x_i \mapsto y_i$ of the local ring $R:=\Q[[x_1,x_2,\ldots,x_n]]$ is given by $n$-tuples of series
$(y_1,y_2,\ldots,y_n)$ with the property that
\[\left| \frac{\partial y}{\partial x}\right| =\left| \begin{array}{ccc}\frac{\partial y_1}{\partial x_1} &\ldots&\frac{\partial y_1}{\partial x_n}\\
\ldots&\ldots&\ldots\\
\frac{\partial y_n}{\partial x_1} &\ldots&\frac{\partial y_n}{\partial x_n}\\
\end{array} \right| \notin \cM:=(x_1,x_2,\ldots,x_n) \subset R\]

One has the following {\em Formal Morse lemma}: there exist an automorphism $\phi$ of $R$ 
such that
\[ \phi(f)=f_2\]
Such a $\phi$ is obtained by an iteration: if \[f=f_2+f_k+f_{k+1}+\ldots,\]
then we can find an automorphism $\phi_k$ such that 
\[\phi_k(f)=f_2+\tilde{f}_{k+1}+\ldots\]
To find $\phi_k$ it is sufficient to write $f_k=\sum a_i \partial f/\partial x_i$ and set
$\phi_k(x_i)=x_i-a_i$.\\

Alternatively, we may say that one can find formal coordinates
$y_i=\phi(x_i)$ such that 
\[f_2+f_3+\ldots=y_1^2+y_2^2+\ldots+y_n^2\]

By the transformation formula for integrals we get
\[ \int_{\Gamma(t)} \omega=\int_{B(t)} \phi^*(\omega)\]
When we write
\[\omega:=A(x)dx_1dx_2\ldots dx_n\]
then 
\[ \phi^*(\omega)=A(x(y)) \left| \frac{\partial x(y)}{\partial y} \right| dy_1dy_2\ldots dy_n\]
which can be expanded in a series in the coordinates $y_i$ as
\[\phi^*(\omega) = \sum_{\alpha} J_{\alpha} y^{\alpha} dy_1 dy_2 \ldots dy_n\]
where the $J_{\alpha} \in \Q$. So we get
\[\int_{\Gamma(t)} \omega=\int_{B(t)}\phi^*(\omega) = \sum_{\alpha} J_{\alpha} \int_{B(t)}y^{\alpha} dy_1 dy_2 \ldots dy_n\]
The integrals
\[I({\alpha}):=\int_{B(t)} y^{\alpha} dy_1dy_2\dots dy_n\]
can be reduced to the generalised Beta-integral and one has:

\begin{lemma}
\rule{0mm}{1mm}\\
(i) 
\[I(\alpha_1,\alpha_2,\ldots,\alpha_n) =0\]
when some $\alpha_i$ is odd.\\
(ii)
\[I(2k_1,2k_2,\ldots,2k_n)=\frac{\Gamma(k_1+1/2)\Gamma(k_2+1/2)\ldots \Gamma(k_n+1/2)}{\Gamma(k_1+k_2+\ldots+n/2+1)} t^{k_1+k_2\ldots+k_n+n/2}
\]  
\end{lemma}

As a consequence we have

\begin{eqnarray*}
\int_{\Gamma(t)}&=&\sum J_{\alpha} I(\alpha) t^{k_1+k_2+\ldots+k_n+n/2}\\
&=&I(0)t^{n/2}(1+a_1t+a_2t^2+\ldots)
\end{eqnarray*}

The coefficient 
\[I(0) =\frac{\Gamma(1/2)^n}{\Gamma(n/2+1)}\]
is the volume of the $n$-dimensional unit ball. 
As $I(\alpha)/I(0,0,\ldots,0) \in \Q$, we see that the $a_i$ are also
in $\Q$.

So we see that the period integral
 
\[\phi(t) =\frac{d}{dt} \int_{\Gamma(t)} \omega\]

has, up to a prefactor, a series expansion with rational coefficients, that
can be computed algebraically be a very simple, although memory consuming 
algorithm. Pavel Metelitsyn is currently working on an implemenatation.\\

{\bf Aknowledgment:} We would like to thank the organisers for
inviting us to the {\em Workshop on Arithmetic and Geometry of K3
  surfaces and Calabi-Yau threefolds} held in the period $16-25$
August $2011$ at the Fields Institute.   
We also thank M. Bogner and J. Hofmann for help  with analysis of the
examples. Furthermore, I thank  G. Almkvist, W. Zudilin for continued
interest in this  
crazy project. Part of this research was done during the stay of the
first named author  as a
guestprofessor at the \emph{Schwerpunkt Polen} of the Johannes
Gutenberg--Universit\"at in Mainz.

\end{document}